\documentclass[10pt,final]{article}

\usepackage{amssymb,amsmath,amsthm,authblk,graphicx,subfig,url,mathptmx}

\newtheoremstyle{theorem}
  {10pt}		  
  {10pt}  
  {\sl}  
  {\parindent}     
  {\bf}  
  {. }    
  { }    
  {}     
\theoremstyle{theorem}
\newtheorem{theorem}{Theorem}

\newtheoremstyle{defi}
  {10pt}		  
  {10pt}  
  {\rm}  
  {\parindent}     
  {\bf}  
  {. }    
  { }    
  {}     
\theoremstyle{defi}

\begin{document}
	\begin{center}
		{\Large \bf{Curves of Constant Breadth According to Darboux Frame \ }}
	\end{center}
	\centerline{\large B\"{u}lent Altunkaya $^{1}$, Ferdag KAHRAMAN AKSOYAK $^{2}${\footnotetext{
				{ E-mail: $^{1}$bulent\_altunkaya@hotmail.com (B. Altunkaya ); $^{2}$ferda.kahraman@yahoo.com (F. Kahraman Aksoyak)}} }}
	
	\
	\centerline{\it $^{1}$ Ahi Evran University, Division of Elementary Mathematics Education
		K\i r\c{s}ehir, TURKEY}
	
	\centerline{\it $^{2}$ Ahi Evran University, Division of Elementary Mathematics Education
		K\i r\c{s}ehir, TURKEY}

	\begin{abstract}
		In this paper, we investigate constant breadth curves on a surface according
		to Darboux frame and give some characterizations of these curves.
	\end{abstract}
	
	\begin{quote}\small
		{\it{Key words and Phrases}: Darboux frame, constant breadth curve, Euclidean space.}
	\end{quote}
	\begin{quote}\small
		2010 \textit{Mathematics Subject Classification}: 53B25 ; 53C40  .
	\end{quote}

\section{Introduction}

Since the first introduction of constant breadth curves in the plane by L.
Euler in 1778 \cite{euler}, many researchers focused on this subject and found
out a lot of  properties about constant breadth curves in the
plane \cite{mellish}, \cite{ball}, \cite{struik}. Fujiwara has introduced
constant breadth curves, by taking a closed curve whose normal plane at a
point $P$ has only one more point $Q$ in common with the curve and for which
the distance $d\left( P,Q\right) $ is constant.

After the development of cam design, researchers have shown strong
interest to this subject again and many interesting properties have been
discovered. For example K\"{o}se has defined a new concept called space
is a pair of curve of constant breadth in \cite{kose2}, a pair of unit speed space
curves of class $C^{3}$ with non-vanishing curvature and torsion in $E^{3},$
which have parallel tangents in opposite directions at corresponding points,
and the distance between these points is always constant by using the Frenet
frame.

The characterizations of K\"{o}se's paper on constant breadth curves in the
space has led us to investigate this topic according to Darboux frame on a
surface.

\section{Basic Concepts}

Now, we introduce some basic concepts about our study. Let $M$ be an
oriented surface and $\beta $ be a unit speed curve of class $C^{3}$\ on $M.$
As we know, $\beta $ has a natural frame called Frenet frame $\{T,N,B\}$
with properties below:%
\begin{equation}
\begin{array}{rl}
T^{\prime }= & \varkappa N, \\ 
N^{\prime }= & -\varkappa T+\tau B \\ 
B^{\prime }= & -\tau N%
\end{array}%
\end{equation}%
where $\varkappa $ is the curvature, $\tau $ is the torsion, $T$ is the unit
tangent vector field, $N$ is the principal normal vector field and $B$ is
the binormal vector field of the curve $\beta .$

Let us take the unit tangent vector field of the curve $\beta $ and the unit
normal vector field $n$ of the surface $M.$ If we define unit vector field $%
g $ as $g=\left( n\circ \beta \right) \times T$ where $\times $ cross
product. We will have a new frame called Darboux frame $\{T,g,n\circ \beta
\}.$ The relations between these two frames can be given as follows:%
\begin{equation}
\left[ 
\begin{array}{c}
T \\ 
g \\ 
n\circ \beta%
\end{array}%
\right] =\left[ 
\begin{array}{ccc}
1 & 0 & 0 \\ 
0 & \cos \alpha & \sin \alpha \\ 
0 & -\sin \alpha & \cos \alpha%
\end{array}%
\right] \left[ 
\begin{array}{c}
T \\ 
N \\ 
B%
\end{array}%
\right]
\end{equation}%
where $\alpha (s)$ is the angle between the vector fields $n\circ \beta $ and 
$B.$ If we take the derivatives of $T,$ $g,$ $n$ with respect to $s$, we
will have%
\begin{equation}
\left[ 
\begin{array}{c}
T^{\prime } \\ 
g^{\prime } \\ 
\left( n\circ \beta \right) ^{\prime }%
\end{array}%
\right] =\left[ 
\begin{array}{ccc}
0 & k_{g} & k_{n} \\ 
-k_{g} & 0 & t_{g} \\ 
-k_{n} & -t_{g} & 0%
\end{array}%
\right] \left[ 
\begin{array}{c}
T \\ 
g \\ 
n\circ \beta%
\end{array}%
\right]
\end{equation}%
where $k_{g},$ $k_{n}$ and $t_{g}$ are called the geodesic curvature, the
normal curvature and the geodesic torsion respectively. Then, we will have
following relations 

\begin{eqnarray}
k_{g} &=&\varkappa \cos \alpha \\
k_{n} &=&\varkappa \sin \alpha  \notag \\
t_{g} &=&\tau -\alpha ^{\prime }  \notag
\end{eqnarray}%
In the differential geometry of surfaces, for a curve $\beta (s)$ lying on a
surface, there are following cases:

$i)$ $\beta $ is a geodesic curve if and only if $k_{g}=0.$

$ii)$ $\beta $ is an asymptotic line if and only if $k_{n}=0.$

$iii)$ $\beta $ is a principal line if and only if $t_{g}=0.$

\section{Curves of Constant Breadth According to Darboux Frame}

Let $\beta (s)$ and $\beta ^{\ast }(s^{\ast })$ be a pair of unit speed
curves of class $C^{3}$ with non-vanishing curvature and torsion in $E^{3}$
which have parallel tangents in opposite directions at corresponding points
and the distance between these points is always constant. We will call $\left( \beta ^{\ast },\beta \right) $ as curve pair of constant breath.

If $\beta $ lies on a surface, it has Darboux frame in addition to Frenet
frame with properties (1), (2), (3) and (4). So we may write for $\beta
^{\ast }$%
\begin{equation*}
\beta ^{\ast }(s^{\ast })=\beta (s)+m_{1}(s)T(s)+m_{2}(s)g(s)+m_{3}(s)\left(
n\circ \beta \right) (s)
\end{equation*}%
If we differentiate this equation with respect to $s$ and use (3), we will
have%
\begin{eqnarray}
\left( \beta ^{\ast }\right) ^{\prime } &=&\frac{d\beta ^{\ast }}{ds^{\ast }}%
\frac{ds^{\ast }}{ds}=\left( 1+m_{1}^{\prime }-m_{3}k_{n}-m_{2}k_{g}\right)
T+\left( m_{1}k_{g}+m_{2}^{\prime }-m_{3}t_{g}\right) g \\
&&+\left( m_{1}k_{n}+m_{2}t_{g}+m_{3}^{\prime }\right) \left( n\circ \beta
\right)   \notag
\end{eqnarray}%
and 
\begin{equation*}
\frac{d\beta ^{\ast }}{ds}=\frac{d\beta ^{\ast }}{ds^{\ast }}\frac{ds^{\ast }%
}{ds}=T^{\ast }\frac{ds^{\ast }}{ds}
\end{equation*}%
As we know $\left \langle T,T^{\ast }\right \rangle =-1.$ Then 
\begin{equation*}
-\frac{ds^{\ast }}{ds}=1+m_{1}^{\prime }-m_{3}k_{n}-m_{2}k_{g}.
\end{equation*}%
So we find from (5),%
\begin{eqnarray}
m_{1}^{\prime } &=&m_{2}k_{g}+m_{3}k_{n}-1-\frac{ds^{\ast }}{ds} \\
m_{2}^{\prime } &=&m_{3}t_{g}-m_{1}k_{g}  \notag \\
m_{3}^{\prime } &=&-m_{1}k_{n}-m_{2}t_{g}.  \notag
\end{eqnarray}%
Let us denote the angle between the tangents at the points $\beta (s)$ and $%
\beta (s+\bigtriangleup s)$ with $\bigtriangleup \theta .$ If we denote the
vector $T(s+\bigtriangleup s)-T(s)$ with $\bigtriangleup T,$ we know $%
\lim_{\bigtriangleup s\rightarrow 0}\frac{\bigtriangleup T}{\bigtriangleup s}%
=\lim_{\bigtriangleup s\rightarrow 0}\frac{\bigtriangleup \theta }{%
\bigtriangleup s}=\frac{d\theta }{ds}=\varkappa .$ We called the angle of
contingency to the angle $\bigtriangleup \theta $ \cite{struik}. Let us denote
the differentiation with respect to $\theta $ with $"\cdot ".$ By using the
equation $\frac{d\theta }{ds}=\varkappa $, we can writre (6) as follows: 
\begin{eqnarray}
\dot{m}_{1} &=&\rho \left( m_{2}k_{g}+m_{3}k_{n}\right) -f(\theta ) \\
\dot{m}_{2} &=&\rho \left( m_{3}t_{g}-m_{1}k_{g}\right)   \notag \\
\dot{m}_{3} &=&\rho \left( -m_{1}k_{n}-m_{2}t_{g}\right) .  \notag
\end{eqnarray}%
where $\rho =\frac{1}{\varkappa }$, $\rho ^{\ast }=\frac{1}{\varkappa ^{\ast
}}$\ and $\rho +\rho ^{\ast }=f(\theta ).$

Now we investigate curves of constant breadth according to Darboux frame for
some special cases:

\subsection{Case (For geodesic curves)}

Let $\beta $ be non straight line geodesic curve on a surface. Then $%
k_{g}=\varkappa \cos \alpha =0$ and $\varkappa \neq 0$, we get $\cos \alpha
=0$. So it implies that $k_{n}=\varkappa ,$ $t_{g}=\tau .$ By using (7), we
have following differential equation system%
\begin{eqnarray}
\dot{m}_{1} &=&m_{3}-f(\theta ) \\
\dot{m}_{2} &=&m_{3}\varphi   \notag \\
\dot{m}_{3} &=&-m_{1}-m_{2}\varphi   \notag
\end{eqnarray}%
where $\varphi =\frac{\tau }{\varkappa }.$ By using (8), we obtain a
differential equation as follows:%
\begin{equation}
\left( \dddot{m}_{1}+\ddot{f}\right) -\frac{d\varphi }{d\theta }\frac{1}{%
\varphi }\left( \ddot{m}_{1}+m_{1}+\dot{f}\right) +\left( 1+\varphi
^{2}\right) \dot{m}_{1}+\varphi ^{2}f=0
\end{equation}%
We assume that $\left( \beta ^{\ast },\beta \right) $is a curve pair of
constant breadth , then 
\begin{equation*}
\left \Vert \beta ^{\ast }-\beta \right \Vert
^{2}=m_{1}^{2}+m_{2}^{2}+m_{3}^{2}=\text{constant}
\end{equation*}%
which imlplies that 
\begin{equation}
m_{1}\dot{m}_{1}+m_{2}\dot{m}_{2}+m_{3}\dot{m}_{3}=0
\end{equation}%
By combining (8) and (10) then we get 
\begin{equation*}
m_{1}f\left( \theta \right) =0
\end{equation*}

\subsubsection{Case $f\left( \protect \theta \right) =0.$}

We assume that $f\left( \theta \right) =0.$ By using (9), we get%
\begin{equation}
\dddot{m}_{1}-\frac{d\varphi }{d\theta }\frac{1}{\varphi }\left( \ddot{m}%
_{1}+m_{1}\right) +\left( 1+\varphi ^{2}\right) \dot{m}_{1}=0
\end{equation}%
If $\beta $ is a helix curve then $\varphi =\varphi _{0}=$constant. From
(11), we have%
\begin{equation*}
\dddot{m}_{1}+\left( 1+\varphi _{0}^{2}\right) \dot{m}_{1}=0
\end{equation*}%
whose solution is 
\begin{equation*}
m_{1}=\frac{1}{\sqrt{1+\varphi _{0}^{2}}}\left( c_{1}\sin \left( \left(
1+\varphi _{0}^{2}\right) \theta \right) -c_{2}\cos \left( \left( 1+\varphi
_{0}^{2}\right) \theta \right) \right) 
\end{equation*}%
So we can find as $m_{2}=-\frac{1}{\varphi _{0}}\left( m_{1}+\ddot{m}%
_{1}\right) $ and $m_{3}=\dot{m}_{1}$.

\subsubsection{Case $m_{1}=0.$}

We assume that $m_{1}=0,$ then by using (9), we get%
\begin{equation}
\ddot{f}-\frac{d\varphi }{d\theta }\frac{1}{\varphi }\dot{f}+\varphi ^{2}f=0
\end{equation}%
If $\beta $ is a helix curve, then $\varphi =\varphi _{0}=$constant. From
(12) we obtain 
\begin{equation*}
\ddot{f}+\varphi _{0}^{2}f=0
\end{equation*}%
whose solution is 
\begin{equation*}
f(\theta )=c_{1}\cos \left( \varphi _{0}\theta \right) +c_{2}\sin \left(
\varphi _{0}\theta \right) 
\end{equation*}%
Since $m_{1}=0$ it implies that%
\begin{equation*}
m_{3}=f(\theta )
\end{equation*}%
\begin{equation*}
m_{2}=-\frac{\dot{m}_{3}}{\varphi _{0}}
\end{equation*}

\begin{theorem}
Let $\beta $ be a geodesic curve and a helix curve. Let $\left( \beta ,\beta
^{\ast }\right) $ be a pair of constant breadth curve. In that case $\beta
^{\ast }$ can be expressed as one of the following cases:

$i)$ 
\begin{equation*}
\beta ^{\ast }(s^{\ast })=\beta (s)+m_{1}(s)T(s)-\frac{1}{\varphi _{0}}%
\left( m_{1}\left( s\right) +\ddot{m}_{1}\left( s\right) \right) g(s)+\dot{m}%
_{1}\left( s\right) \left( n\circ \beta \right) (s)
\end{equation*}%
where $m_{1}=\frac{1}{\sqrt{1+\varphi _{0}^{2}}}\left( c_{1}\sin \left(
\left( 1+\varphi _{0}^{2}\right) \theta \right) -c_{2}\cos \left( \left(
1+\varphi _{0}^{2}\right) \theta \right) \right) $.

$ii)$%
\begin{equation*}
\beta ^{\ast }(s^{\ast })=\beta (s)-\frac{\dot{f}(\theta )}{\varphi _{0}}%
g(s)+f(\theta )\left( n\circ \beta \right) (s)
\end{equation*}%
where $f(\theta )=c_{1}\cos \left( \varphi _{0}\theta \right) +c_{2}\sin
\left( \varphi _{0}\theta \right) .$
\end{theorem}

\subsection{Case (For asymptotic lines)}

Let $\beta $ be non straight line asymptotic line on a surface. Then $%
k_{n}=\varkappa \sin \alpha =0$ and $\varkappa \neq 0$, we have $\sin \alpha
=0$. So we get $k_{g}=\varepsilon \varkappa ,$ $t_{g}=\tau ,$ where $%
\varepsilon =\pm 1.$ By using (7), we have following differential equation
system%
\begin{eqnarray}
\dot{m}_{1} &=&\varepsilon m_{2}-f(\theta ) \\
\dot{m}_{2} &=&m_{3}\varphi -\varepsilon m_{1}  \notag \\
\dot{m}_{3} &=&-m_{2}\varphi   \notag
\end{eqnarray}%
where $\varphi =\frac{\tau }{\varkappa }.$ By using (13), we obtain a
differential equation as follows:%
\begin{equation}
\left( \dddot{m}_{1}+\ddot{f}\right) -\frac{d\varphi }{d\theta }\frac{1}{%
\varphi }\left( \ddot{m}_{1}+m_{1}+\dot{f}\right) +\left( 1+\varphi
^{2}\right) \dot{m}_{1}+\varphi ^{2}f=0
\end{equation}%
We assume that $\left( \beta ^{\ast },\beta \right) $ is a curve pair of
constant breadth then 
\begin{equation*}
\left \Vert \beta ^{\ast }-\beta \right \Vert
^{2}=m_{1}^{2}+m_{2}^{2}+m_{3}^{2}=\text{constant}
\end{equation*}%
which imlplies that 
\begin{equation}
m_{1}\dot{m}_{1}+m_{2}\dot{m}_{2}+m_{3}\dot{m}_{3}=0
\end{equation}%
By combining (13) and (15) then we get 
\begin{equation*}
m_{1}f\left( \theta \right) =0
\end{equation*}

\subsubsection{Case $f\left( \protect \theta \right) =0$}

We assume that $f\left( \theta \right) =0.$ By using (14), we get%
\begin{equation}
\dddot{m}_{1}-\frac{d\varphi }{d\theta }\frac{1}{\varphi }\left( \ddot{m}%
_{1}+m_{1}\right) +\left( 1+\varphi ^{2}\right) \dot{m}_{1}=0
\end{equation}%
If $\beta $ is a helix curve then $\varphi =\varphi _{0}=$constant. From
(16), we have%
\begin{equation*}
\dddot{m}_{1}+\left( 1+\varphi _{0}^{2}\right) \dot{m}_{1}=0
\end{equation*}%
whose solution is 
\begin{equation*}
m_{1}=\frac{1}{\sqrt{1+\varphi _{0}^{2}}}\left( c_{1}\sin \left( \left(
1+\varphi _{0}^{2}\right) \theta \right) -c_{2}\cos \left( \left( 1+\varphi
_{0}^{2}\right) \theta \right) \right) 
\end{equation*}%
So we can find as $m_{2}=\varepsilon \dot{m}_{1}$ and $m_{3}=\varepsilon 
\frac{1}{\varphi _{0}}\left( m_{1}+\ddot{m}_{1}\right) $.

\subsubsection{Case $m_{1}=0$}

We assume that $m_{1}=0.$Then by using (14), we get%
\begin{equation}
\ddot{f}-\frac{d\varphi }{d\theta }\frac{1}{\varphi }\dot{f}+\varphi ^{2}f=0
\end{equation}%
If $\beta $ is a helix curve, then $\varphi =\varphi _{0}=$constant. From
(17) we obtain 
\begin{equation*}
\ddot{f}+\varphi _{0}^{2}f=0
\end{equation*}%
whose solution is 
\begin{equation*}
f(\theta )=c_{1}\cos \left( \varphi _{0}\theta \right) +c_{2}\sin \left(
\varphi _{0}\theta \right) 
\end{equation*}%
Since $m_{1}=0$ it implies that%
\begin{equation*}
m_{2}=\varepsilon f(\theta )
\end{equation*}%
\begin{equation*}
m_{3}=\frac{\dot{m}_{2}}{\varphi _{0}}
\end{equation*}%
where $\theta =\int \varkappa ds.$

\begin{theorem}
Let $\beta $ be an asymptotic line and a helix curve. Let $\left( \beta
,\beta ^{\ast }\right) $ be a pair of constant breadth curve. In that case $%
\beta ^{\ast }$ can be expressed as one of the following cases:

$i)$ 
\begin{equation*}
\beta ^{\ast }(s^{\ast })=\beta (s)+m_{1}(s)T(s)+\varepsilon \dot{m}%
_{1}\left( s\right) g(s)-\frac{\varepsilon }{\varphi _{0}}\left( m_{1}\left(
s\right) +\ddot{m}_{1}\left( s\right) \right) \left( n\circ \beta \right) (s)
\end{equation*}%
where $m_{1}=\frac{1}{\sqrt{1+\varphi _{0}^{2}}}\left( c_{1}\sin \left(
\left( 1+\varphi _{0}^{2}\right) \theta \right) -c_{2}\cos \left( \left(
1+\varphi _{0}^{2}\right) \theta \right) \right) $.

$ii)$%
\begin{equation*}
\beta ^{\ast }(s^{\ast })=\beta (s)+\varepsilon f(\theta )g(s)+\varepsilon 
\frac{\dot{f}(\theta )}{\varphi _{0}}\left( n\circ \beta \right) (s)
\end{equation*}%
where $f(\theta )=c_{1}\cos \left( \varphi _{0}\theta \right) +c_{2}\sin
\left( \varphi _{0}\theta \right) .$
\end{theorem}

\subsection{Case (For principal line)}

We assume that $\beta $ is a principal line. Then we have $t_{g}=0$ and it
implies that $\tau =\alpha ^{\prime }$. By using (7), we get%
\begin{eqnarray}
\dot{m}_{1} &=&m_{2}\cos \alpha +m_{3}\sin \alpha -f(\theta ) \\
\dot{m}_{2} &=&-m_{1}\cos \alpha   \notag \\
\dot{m}_{3} &=&-m_{1}\sin \alpha   \notag
\end{eqnarray}%
By using (18), we obtain following differential equation%
\begin{equation}
\left( \dddot{m}_{1}+\dot{m}_{1}\right) +\left( \dot{m}_{1}+f\right) \dot{%
\alpha}^{2}-\left( \sin \alpha \int m_{1}\cos \alpha d\theta -\cos \alpha
\int m_{1}\sin \alpha d\theta \right) \ddot{\alpha}+\ddot{f}=0
\end{equation}%
since $t_{g}=0$ we obtain $\dot{\alpha}=\frac{\tau }{\varkappa }$. \ We
assume that $\left( \beta ^{\ast },\beta \right) $ is a pair of curve is
constant breadth. In that case 
\begin{equation*}
\left \Vert \beta ^{\ast }-\beta \right \Vert
^{2}=m_{1}^{2}+m_{2}^{2}+m_{3}^{2}=\text{constant}
\end{equation*}%
which imlplies that 
\begin{equation}
m_{1}\dot{m}_{1}+m_{2}\dot{m}_{2}+m_{3}\dot{m}_{3}=0
\end{equation}%
By combining (18) and (20) then we get 
\begin{equation*}
m_{1}f\left( \theta \right) =0
\end{equation*}

\subsubsection{Case $f(\protect \theta )=0$}

We assume that $f\left( \theta \right) =0.$ By using (19), we get%
\begin{equation}
\left( \dddot{m}_{1}+\dot{m}_{1}\right) +\dot{m}_{1}\dot{\alpha}^{2}-\left(
\sin \alpha \int m_{1}\cos \alpha d\theta -\cos \alpha \int m_{1}\sin \alpha
d\theta \right) \ddot{\alpha}=0
\end{equation}%
If $\beta $ is a helix curve then $\dot{\alpha}=\frac{\tau }{\varkappa }=$%
constant. From (21), we have%
\begin{equation*}
\dddot{m}_{1}+\left( 1+\dot{\alpha}^{2}\right) \dot{m}_{1}=0
\end{equation*}%
Then we get 
\begin{equation*}
m_{1}(s)=\frac{1}{\sqrt{1+\dot{\alpha}^{2}}}\left( c_{1}\sin \left( \left( 1+%
\dot{\alpha}^{2}\right) \theta \right) -c_{2}\cos \left( \left( 1+\dot{\alpha%
}^{2}\right) \theta \right) \right) 
\end{equation*}%
By using (18) we obtain 
\begin{eqnarray*}
\dot{m}_{2} &=&-m_{1}\cos \alpha  \\
\dot{m}_{3} &=&-m_{1}\sin \alpha 
\end{eqnarray*}%
where $\alpha =\int \tau ds.$

\begin{theorem}
Let $\beta $ be a principal line and a helix curve. Let $\left( \beta ,\beta
^{\ast }\right) $ be a pair of constant breadth curve such that $%
\left \langle \beta ^{\ast },T\right \rangle =m_{1}\neq 0.$ In that case $%
\beta ^{\ast }$ can be expressed as:%
\begin{equation*}
\beta ^{\ast }=\beta +m_{1}(s)T(s)-m_{1}(s)\cos \alpha g(s)-m_{1}\sin \alpha
\left( n\circ \beta \right) (s)
\end{equation*}%
where $m_{1}(s)=\frac{1}{\sqrt{1+\dot{\alpha}^{2}}}\left( c_{1}\sin \left(
\left( 1+\dot{\alpha}^{2}\right) \theta \right) -c_{2}\cos \left( \left( 1+%
\dot{\alpha}^{2}\right) \theta \right) \right) .$
\end{theorem}

\subsubsection{Case $m_{1}=0$}

If $m_{1}=0,$ then from (19) 
\begin{equation}
\ddot{f}+\dot{\alpha}^{2}f=0
\end{equation}%
where $\dot{\alpha}=\frac{\tau }{\varkappa }.$On the other hand since $%
m_{1}=0$ from (18) we have $m_{2}=c_{2}$=constant, $m_{3}=c_{3}$=constant
and 
\begin{equation}
f=c_{2}\cos \alpha +c_{3}\sin \alpha 
\end{equation}%
By combining (22) and (23)%
\begin{equation*}
\ddot{\alpha}\left( -c_{2}\sin \alpha +c_{3}\cos \alpha \right) =0
\end{equation*}%
In that case, if $\ddot{\alpha}=0$ then we obtain that $\dot{\alpha}=\frac{%
\tau }{\varkappa }=$constant. $\beta $ becomes a helix curve. If $-c_{2}\sin
\alpha +c_{3}\cos \alpha =0$ then we have $\alpha =$constant. This means
that $\beta $ is a planar curve.

\begin{theorem}
Let $\beta $ be a principal line. Let $\left( \beta ,\beta ^{\ast }\right) $
be a pair of constant breadth curve such that $\left \langle \beta ^{\ast
},T\right \rangle =m_{1}=0.$ In that case $\beta $ is a helix curve or a
planar curve and $\beta ^{\ast }$ can be expressed as:%
\begin{equation*}
\beta ^{\ast }=\beta +c_{2}g(s)+c_{3}\left( n\circ \beta \right) (s)
\end{equation*}
\end{theorem}


\begin{thebibliography}{99}
\bibitem{koca2} Kocayi\u{g}it, H., Zennure, \c{C}., Same characterizations
of of constant breadth spacelike curves in Minkowski 4-space $E_{1}^{4}$,
NTMSCI 3 (2)
\end{thebibliography}
\end{document}